\documentclass[10pt]{article}
\usepackage{amsmath}
\usepackage{graphicx} 
\usepackage{amssymb} 
\usepackage{amsthm}
\usepackage{subcaption} 
\usepackage{tikz}
\usepackage{authblk}
\usepackage{hyperref}
\usepackage{xcolor}

\hypersetup{
    colorlinks = true,
    allcolors=.
}
\usetikzlibrary{shapes}
\usetikzlibrary{positioning}
\usetikzlibrary{arrows}
\usetikzlibrary{decorations.markings}

\theoremstyle{plain}

\newtheorem{theorem}{Theorem}[section]
\newtheorem{lemma}[theorem]{Lemma}
\newtheorem{corollary}[theorem]{Corollary}

\DeclareMathOperator{\tr}{tr}

\newcommand{\pa}[1]{\left(#1\right)} 
\newcommand{\mat}[1]{\begin{matrix}#1\end{matrix}} 
\newcommand{\pmat}[1]{\pa{\mat{#1}}} 

\title{Hermitian adjacency matrices with at most three distinct eigenvalues}
\author[1]{Saieed Akbari\thanks{The research visit of S. Akbari at Simon Fraser University was supported in part by the ERC Synergy grant (European Union, ERC, KARST, project number 101071836).}}
\author[2]{Jonathan Aloni}
\author[2]{Maxwell Levit\thanks{Postdoctoral fellowship at SFU supported through NSERC Discovery Grants R832714 and R611368 (Canada), and the ERC Synergy grant KARST (European Union, ERC, KARST, project number 101071836).}}
\author[2]{Bojan Mohar\thanks{Supported in part by the NSERC Discovery Grant R832714 (Canada), and in part by the ERC Synergy grant KARST (European Union, ERC, KARST, project number 101071836).
On leave from FMF, Department of Mathematics, University of Ljubljana.}}
\author[2]{Steven Xia\thanks{Supported by the NSERC Undergraduate Student Research Award 615898.}}
\affil[1]{Sharif University of Technology}
\affil[2]{Simon Fraser University}

\date{}

\begin{document}

\maketitle

\begin{abstract}
    We study oriented graphs whose Hermitian adjacency matrices of the second kind have few eigenvalues. We give a complete characterization of the oriented graphs with two distinct eigenvalues, showing that there are only four such graphs. We extend this result to mixed graphs. We show that there are infinitely many regular tournaments with three distinct eigenvalues. We extend our main results to Hermitian adjacency matrices defined over other roots of unity. 
\end{abstract} 

\section{Introduction}
Let $D$ be an oriented graph and let $\omega=\frac{1+i\sqrt{3}}{2}$, a primitive $6$th root of unity. The \textbf{Hermitian adjacency matrix} $H_\omega(D)$ is defined by 

\[H_\omega(D)_{u,v}=\begin{cases}
    \omega&   u\rightarrow v \\
    \overline{\omega}&  u\leftarrow v \\
    0&\mbox{otherwise}.
\end{cases} \]

Our first main result is a complete characterization of the oriented graphs $D$ for which $H_\omega(D)$ has exactly two distinct eigenvalues. In particular, we show that there are only four such graphs.

\begin{figure}
    \centering
    \includegraphics[scale=0.66]{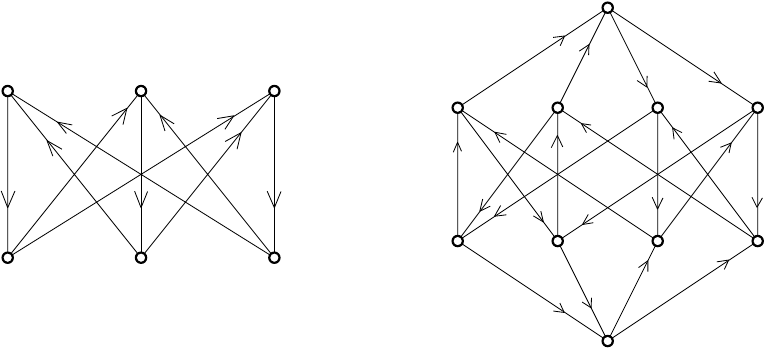}
    \caption{The special oriented graphs from Theorem \ref{Thm:2ev_Char} (c) and (d).}
    \label{fig:Herm_mats}
\end{figure}

\begin{theorem}\label{Thm:2ev_Char}
If $D$ is a connected oriented graph and $H_\omega(D)$ has exactly two distinct eigenvalues, then $D$ is one of the following:
\begin{itemize}
    \item[(a)] The directed edge with eigenvalues $\{1,-1\}$.
    \item[(b)] The directed triangle with eigenvalues $\{1,-2\}$.
    \item[(c)] A special orientation of $K_{3,3}$ with eigenvalues $\{\sqrt{3},-\sqrt{3}\}$; see Figure~\ref{fig:Herm_mats}.
    \item[(d)] A special orientation of $K_{5,5}$ minus a perfect matching with eigenvalues $\{2,-2\}$; see Figure~\ref{fig:Herm_mats}.
\end{itemize}
\end{theorem}

This was not the result we expected at the outset, given that analogous problems for signed graphs admit numerous infinite families and have so far evaded full characterizations. We refer the reader to the papers \cite{G17,HTW19,R22,ST91,S22} for various infinite families of signed graphs with two distinct eigenvalues.

In the remainder of this section we give some context for the problem and the required preliminary definitions. In Section \ref{Sec:two_ev} we prove Theorem \ref{Thm:2ev_Char}. We then extend the characterization to mixed graphs with two eigenvalues obtaining only $C_4$ and trivial mixed orientations of $K_n$ as additional examples. In Section  \ref{Sec:3ev} we show that, in counterpoint to Theorem \ref{Thm:2ev_Char}, there are infinitely many directed graphs with exactly three distinct $H_\omega$-eigenvalues. In particular, we show that there are infinitely many regular tournaments with this property by generalizing results of \cite{dCGKPM92} which leverage skew-symmetric Hadamard matrices. In Section 4 we discuss the extension of Section \ref{Sec:two_ev} to Hermitian adjacency matrices $H_\sigma$ for $\sigma$ some primitive $k$th root of unity. We show that the case $k=3$ is essentially identical to $k=6$, while the case $k=4$ is closely related to signed graphs with two eigenvalues. We conclude with a discussion for general $k$ and prove that for $k>8$ the directed edge is the only oriented graph with two $H_\sigma$-eigenvalues when $\sigma=\cos(\frac{2\pi}{k})+i\sin(\frac{2\pi}{k})$.

\subsection{Context for the Problem}

Hermitian adjacency matrices for digraphs were first studied by Guo and Mohar \cite{GM15} and independently by Liu and Li \cite{LL15}. In both of these papers, the authors used the primitive fourth root of unity $i$.

Subsequently, Mohar \cite{M20} introduced the variant with $\omega=\frac{1+i\sqrt{3}}{2}$ which we are using at present. More precisely, he introduced a more general definition for digraphs which specializes to our definition in the case of oriented graphs. Subsequently, Li and Yu characterized the mixed graphs with  $H_\omega$-eigenvalues in certain small intervals, \cite{LY21}. Together with Zhou \cite{ZLY23}, they then characterized the mixed graphs with the least $H_\omega$-eigenvalue greater than $-\frac{3}{2}$. 

All of these results belong to the broader study of spectral properties of the complex unit gain graphs defined by Reff \cite{R11}. Gain graphs with two distinct eigenvalues are related to many interesting and exceptional combinatorial objects such as distance-regular covering graphs and systems of equiangular or biangular lines in complex space, see the papers \cite{GLS21} and \cite{vDW22b} for various perspectives on these connections.

\subsection{Why Sixth Roots of Unity?}

The reader may reasonably wonder why we, and others, have chosen to work specifically with  $\omega=\frac{1+i\sqrt{3}}{2}$. First, note that there are only two primitive sixth roots of unity, $\omega$ and $\overline{\omega}$ and that $H_\omega(D)=H^T_{\overline{\omega}}(D).$ So the two choices of primitive sixth root give cospectral matrices and we lose nothing by working specifically with $\omega$. But why sixth roots?

The motivation for this choice is presented in \cite{M20}, where it is pointed out that $\omega+\overline{\omega}=1$ perfectly models the convention for mixed graphs that having two directed arcs in opposite directions is equivalent to an undirected edge. This choice has useful consequences.

For instance, using the determinant expansion of the characteristic polynomial of $H_\omega(D)$ one can verify that the coefficients of this characteristic polynomial are integers for any mixed graph $D$. See Theorem 2.6 from \cite{LY21} for an explicit description of these coefficients. On the other hand, let $D$ be a triangle, directed so that some vertex has in-degree 2, and let $\sigma=\cos(\frac{2\pi}{k})+i\sin(\frac{2\pi}{k})$ be a primitive $k$th root of unity. The constant term of the characteristic polynomial of $H_\sigma(D)$ is $-\sigma-\overline{\sigma}=-2\mbox{Re}(\sigma)=-2\cos(\frac{2\pi}{k})$ which is an integer only if $k\in \{1,2,3,4,6\}$.

In Section \ref{Sec:other_roots} we will highlight some specific instances in our proofs where the choice $\omega=\frac{1+i\sqrt{3}}{2}$ is important.

\subsection{Preliminary Notation and Definitions}

An \textbf{oriented graph} $D=(V,A)$ consists of a vertex set $V$ and a set $A$ of ordered pairs of vertices called the \textbf{arcs}. For each vertex $v\in V$ we require $(v,v)\notin A$ so that $D$ is loopless, and for each pair of vertices $u$ and $v$ at most one of $(u,v)$ or $(v,u)$ is an arc of $A$. We will write $u \rightarrow v$ or $v\leftarrow u$ to denote $(u,v)\in A$. Moreover, we extend the previous notation so that $x\rightarrow y\leftarrow v$ denotes ``both $(x,y)\in A$ and $(v,y)\in A$''. The \textbf{underlying graph} of $D$ is the simple graph $G=(V,E)$ with $E=\{\{u,v\}:(u,v)\in A \mbox{ or }(v,u)\in A\}$. We will occasionally abuse notation and write $|E(D)|$ to mean $|A(D)|$ to avoid confusion with the standard notation of adjacency matrices. In an undirected graph we let $N(v)$ denote the neighbourhood of $v$. We say that two vertices of $D$ are \textbf{adjacent} if they are adjacent in the underlying graph of $D$. We say that an oriented graph $D$ is connected if the underlying graph of $D$ is connected. The \textbf{in-degree} (respectively, \textbf{out-degree}) of a vertex $v$ is the size of the set $\{u: (u,v)\in A\}$ (respectively, $\{u: (v,u)\in A\}$) and is denoted $d^-(v)$ (respectively $d^+(v)$). The \textbf{degree} of $v$ is $d(v)=d^+(v)+d^-(v).$ An oriented graph is \textbf{regular} if all in-degrees and all out-degrees of all vertices are equal. A \textbf{walk} in a graph is a sequence of consecutively adjacent vertices. A \textbf{$k$-walk} in an oriented graph is a sequence $v_1,v_2,\dots,v_{k+1}$ of vertices in which the vertices $v_i$ and $v_{i+1}$ in each successive pair are adjacent. Note that we do not stipulate the direction of the arc connecting $v_i$ to $v_{i+1}$. If we do have the additional feature that each of $(v_i,v_{i+1})\in A$, then we have a \textbf{directed walk}. The \textbf{value} of a walk $u_0u_1\cdots u_r$ in $D$ is $\prod_{i=1}^rH_{u_{i-1}u_i}$. A \textbf{directed cycle} is a regular oriented cycle. 

We recall the interlacing theorem for Hermitian matrices (see, e.g. \cite[p.~246]{HJ12}). We will use it several times.

\begin{theorem}[Eigenvalue interlacing]
         Let $\theta_i(M)$ denote the $i$th largest eigenvalue of a Hermitian matrix $M$. If $B$ is an $m\times m$ induced submatrix of the $n\times n$ Hermitian matrix $A$, then for all $1\leq i\leq m$, the eigenvalues of $A$ and $B$ satisfy
         \[\theta_{i+n-m}(A)\leq \theta_i(B)\leq \theta_i(A).\]
\end{theorem}

\section{Two Distinct Eigenvalues}\label{Sec:two_ev}

In this section, we characterize the oriented graphs with precisely two distinct $H_{\omega}$-eigen\-values. Our main result is Theorem \ref{Thm:2ev_Char}, where we show that there are exactly four such graphs. In Theorem \ref{Thm:Mixed_2ev} we extend the result to mixed graphs.

From now through Section \ref{Sec:rgeq3}, $D$ will be an oriented graph on $n$ vertices for which $H=H_{\omega}(D)$ has two distinct eigenvalues: $r$ with multiplicity $m$, and $s$ with multiplicity $n-m$. We assume $s<0<r$. Our proof of Theorem \ref{Thm:2ev_Char} consists of several small claims, which are organized as follows:

\begin{enumerate}
    \item We first show that $s\geq -2$ with equality if and only if $D$ is regular. We then characterize the irregular examples, for which there are only two possible underlying graphs, $K_2$ and $K_{3,3}$.

    \item Next we consider the regular case. Here, $r$ is equal to the out-degree of each vertex. We show that there is a single example for $r=1$, the directed triangle. We then show
    that the case $r=2$ is equivalent to the underlying graph being triangle-free. We show that $K_{5,5}$ minus a perfect matching is the only possible underlying graph in this case. 

    \item Finally, we give combinatorial arguments which rule out the possibility of $r>2$.  
\end{enumerate}

\subsection{Irregular Oriented Graphs}

First we give a lower bound on the negative eigenvalue $s$.

\begin{lemma} \label{Lem:s_geq2}
    Let $D$ be an oriented graph with two distinct $H_{\omega}$-eigenvalues $r$ and $s$, with $r>s$.
    Then $s\geq-2$, with equality if and only if $D$ is $r$-regular.
\end{lemma}

\begin{proof}
    Set $H=H_\omega(D)$ and $n=|V(D)|$. Let $j$ denote the all-ones vector of length $n$. Since the entries of $H$ are $\omega$ and $\bar{\omega}$ with $\omega+\bar{\omega}=1$ and since $r$ is the largest eigenvalue of $H$ we have
    \begin{equation}\label{eq:rayley} r\geq \frac{j^\intercal Hj}{j^\intercal j}=\frac{1}{n}\sum_{i>j}(H_{i,j}+H_{j,i})=\frac{|E(D)|}{n}\end{equation} with equality if and only if $D$ is regular.

    Now we calculate $\tr(H^2)$ in two ways. On the one hand, the minimal polynomial of $H$ is $(x-r)(x-s)=x^2-(r+s)x+rs$, so 
    \begin{equation}\label{Eq:Min_poly1}
        H^2-(r+s)H+rsI=0.
    \end{equation} 

    Since $H$ has zero diagonal, this implies $\tr(H^2)=-nrs.$ On the other hand, $\tr(H^2)$ is the sum, over all vertices of $D$, of the value of 2-walks from $v$ to $v$. There is precisely one such walk for each arc incident with $v$, and the value of each such walk is $\omega\overline{\omega}=1$, so we obtain \begin{equation} \label{eq:nrs}
        -nrs=\tr(H^2)=\sum_{v\in V} d(v)=2|E(D)|.
    \end{equation}

    Combining Equations (\ref{eq:rayley}) and  (\ref{eq:nrs}) gives the desired bound.
     \end{proof}

The next lemma is a Hermitian variant of Lemma 2.1 from \cite{R22}.
\begin{lemma} \label{Lem:min_poly}
    Let $D$ be an oriented graph with two distinct $H_{\omega}$-eigenvalues, $r$ and $s$.
    Then the underlying graph of $D$ is $|rs|$-regular. Moreover, either $r,s\in\mathbb{Z}$ or $r=-s=\sqrt{k}$ for some $k\in\mathbb{N}$.
\end{lemma}

\begin{proof}
    Let $H=H_\omega(D).$ We again consider the minimal polynomial of $H$ to obtain Equation (\ref{Eq:Min_poly1}). 
    This implies the underlying graph of $D$ is $|rs|$-regular. 
    
    Let $n=|V(D)|$ and let $m$ be the multiplicity of $r$. If $m=\frac{n}2$, then since $\tr(H)=0$ we have $r=-s$ and
    \[H^2-(r+s)H+rsI=H^2+(-s)sI=0.\]
   Thus: $H^2=s^2I$, the underlying graph of $D$ is $s^2$-regular, $s^2\in\mathbb{Z}$, and \[r=-s=\sqrt{s^2}.\]

    Now assume $m\neq\frac{n}2$ and pick vertices $u,v$ so that $u\rightarrow v$. We consider the $(u,v)$ entry of the matrix equation (\ref{Eq:Min_poly1}).
    The $(u,v)$ entry of $H^2$ is the sum of the values of all 2-walks from $u$ to $v$. Since $\omega^4=\omega^{-2}$ we may write this value as $a+b\omega^2+c\omega^4$ for some non-negative integers $a,b,c$.
    
    So we have \[a+b\omega^2+c\omega^4=(r+s)\omega.\] Since $1+\omega^2+\omega^4=0$ we may rewrite the above as \[(b-a)\omega^2+(c-a)\omega^4=(r+s)\omega.\] Multiplying both sides by $\omega^{-1}$ and noting that $\omega^3=-1$ we have \[(b-a)\omega=(a-c)\omega^3+(r+s)=(c-a)+(r+s).\] Since $r$ and $s$ are the eigenvalues of a Hermitian matrix they are real, so the right-hand side is real. The left-hand side is an integer multiple of $\omega$, so we must have $b-a=0$, hence $r+s=a-c$ is an integer.
    
    Now, substituting $s=\frac{mr}{m-n}$ into $r+s$ and factoring out $r$ we have 
    \[r\left(1+\frac{m}{m-n}\right)\in\mathbb{Z}.\] Since $m\neq \frac{n}{2}$ we have  $1+\frac{m}{m-n}\neq 0$, so $r\in\mathbb{Q}$. Since $r$ is an algebraic integer, $r\in\mathbb{Z}$, and thus $s\in\mathbb{Z}$ as well.
\end{proof}

Now we come to a key characterization of the structure of the underlying graph.
    
\begin{lemma} \label{Lem:3common}
If $G$ is the underlying graph of the oriented graph $D$ with exactly two distinct $H_\omega$-eigenvalues $r,s$, where $r=-s$, then for any pair of vertices $u,v$ in $G$, the quantity $|N(u)\cap N(v)|$ is a multiple of $3$.
\end{lemma}

\begin{proof}  
From the minimal polynomial of $H$ we find that  $H^2=s^2I$. Hence for each pair $\{u,v\}$ of distinct vertices of $D$, the sum of the values on 2-walks from $u$ to $v$ is 0. For a pair of distinct vertices, the 2-walks between them take values in $\{1,\omega^2,\omega^4\}$ (with two different orientations that allow for the value 1). Hence there are non-negative integers $a,b,c$ so that $a+b\omega^2+c\omega^4=0$. Since $1+\omega^2+\omega^4=0$, we may rewrite the above as \[(b-a)\omega^2+(c-a)\omega^4=0.\] Multiplying both sides by $\omega^2$ we find \[(b-a)\omega^4=a-c.\] The right-hand side is real, so the left-hand side is zero, so the right-hand side is zero. We conclude that $a=b=c$, so $u$ and $v$ have $3a$ common neighbours. 
\end{proof}

\begin{lemma} \label{Lem:k33} 
If $G$ is a connected $3$-regular graph and each pair of vertices have zero or three common neighbors, then $G$ is isomorphic to $K_{3,3}$. 
\end{lemma}

\begin{proof}
If there is no pair of vertices at distance 2, then $G$ is complete, but since $G$ is 3-regular it must be $K_4$, and any pair of distinct vertices have only two common neighbors, a contradiction. So there is a pair of vertices $\{u,v\}$ at distance 2. They have at least one common neighbour, hence they have exactly three common neighbours $\{x,y,z\}$. Now, $x$ and $y$ have common neighbours $u$ and $v$, so they must share another common neighbour. This third common neighbour cannot be $z$, because then $z$ would have degree at least 4. So there is a sixth vertex $w$ incident with both $x$ and $y$. By the same argument, $x$ and $z$ must have a third common neighbour, and since $x$ already has degree $3$, this neighbour must be $w$. Hence $G$ is a $3$-regular graph containing $K_{3,3}$ as a subgraph. Since $G$ is connected we are done.
\end{proof}

\begin{theorem} \label{Thm:K33}
    Let $D$ be a connected oriented graph with two distinct $H_{\omega}$-eigenvalues $r$ and $s$, where $r>s$.
    If $s>-2$, then $D$ is either a single directed edge or the orientation of $K_{3,3}$ from Figure \ref{fig:Herm_mats}. 
\end{theorem}

\begin{proof}
    Since $s>-2$, Lemma \ref{Lem:min_poly} implies $s\in \{-1,-\sqrt{2},-\sqrt{3}\}.$

    If $s=-1$, then $D$ has at least two vertices. Indeed, the connected oriented graph on two vertices has two distinct eigenvalues $1$ and $-1$.
    However, if $D$ has more than two vertices then it contains an induced subgraph on three vertices. All connected oriented graphs on three vertices have minimum $H_\omega$-eigenvalue less than $-1$, so interlacing implies that no graph on more than two vertices has $s=-1$.

    If $s=-\sqrt{2}$ then Lemma \ref{Lem:min_poly} implies $r=\sqrt{2}$ and the underlying graph of $D$ is $2$-regular, a cycle. But this contradicts Lemma \ref{Lem:3common}.
    
    If $s=-\sqrt{3}$ then Lemma \ref{Lem:min_poly} implies $r=\sqrt{3}$ and the  underlying graph is $3$-regular. Lemma \ref{Lem:3common} implies that each pair of vertices have zero or three common neighbours, so Lemma \ref{Lem:k33} implies the underlying graph is $K_{3,3}$.  
    
    Now, consider a copy of $K_{3,3}$ with one part labeled $u,v,w$ and the other part labeled $x,y,z$. From the proof of Lemma \ref{Lem:3common} we see that among the three paths $(u,x,v)$, $(u,y,v)$ and $(u,z,v)$ exactly one must be directed from $u$ to $v$ and one must be directed from $v$ to $u$. Without loss of generality, assume we have $u\rightarrow x \rightarrow v$ and $v\rightarrow z \rightarrow u.$  The remaining path must give a walk with value 1 so we either have $u\rightarrow y \leftarrow v$ or $u\leftarrow y \rightarrow v$.

    If we have $u\rightarrow y \leftarrow v$ then since two arcs leave $u$ we must have $w\rightarrow z$ in order for there to be a directed 2-walk from $w$ to $u$. Similarly, we must have $w\rightarrow x$ in order for there to be a directed 2-walk from $w$ to $v$. The remaining arc must be $y\rightarrow w$ in order for there to be a directed 2-walk from $u$ to $w$ and we have determined an orientation. 

    If on the other hand we had $u\leftarrow y \rightarrow v$ then a very similar argument implies $x\rightarrow w$, $z\rightarrow w$ and $w\rightarrow y$. The orientation determined in this case is isomorphic to the one in the previous case. 
\end{proof}

\subsection{Regular Oriented Graphs with $r\leq 2$}

It remains to consider connected $r$-regular oriented graphs, the first case to consider is $r=1$, where $D$ must be a directed cycle.

\begin{lemma} A connected $1$-regular oriented graph $D$ with two distinct $H_\omega$-eigenvalues is a directed triangle. \end{lemma}

\begin{proof}
    Certainly $D$ must be a directed cycle. From Lemma \ref{Lem:s_geq2} the eigenvalues are $1$ and $-2$. If $D$ is not a triangle then $D$ contains an induced path on three vertices. Regardless of the orientation, $\sqrt{2}$ is an $H_\omega$-eigenvalue of this induced path, contradicting interlacing. 
\end{proof}

We will see that the case $r=2$ is equivalent to a lack of triangles in the underlying graph. We will also see, somewhat surprisingly, that there is exactly one example in this case. 

\begin{lemma}\label{Lem:r=2}
    Let $D$ be an $r$-regular oriented graph with exactly two distinct $H_{\omega}$-eigenvalues.  If the underlying graph of $D$ is triangle-free, then $r=2$.
\end{lemma}

\begin{proof}
    From Lemma \ref{Lem:s_geq2} we have $s=-2$. The formula for $\tr(D)$ gives 
    \begin{align*}
    0&=\tr(D)=rm-2(n-m).       
    \end{align*}

    The lack of triangles gives 
    \begin{align*}
    0&=\tr(D^3)=r^3m-8(n-m).       
    \end{align*}

    Since $n-m\neq 0$, substituting $rm=2(n-m)$ into the second equation gives $r^2=4$, hence $r=2$.
    \end{proof}
We let $K_{5,5}-M$ denote the graph obtained from $K_{5,5}$ by deleting the edges of a perfect matching.
\begin{lemma} \label{Lem:k55} 
If $G$ is a connected $4$-regular triangle-free undirected graph and each pair of vertices at distance $2$ have three common neighbors, then $G$ is isomorphic to $K_{5,5}-M.$
\end{lemma}

\begin{proof}
Pick a vertex $v$ in $G$ and denote its neighbours $w,x,y,z$. Let $S$ denote the set of vertices at distance 2 from $v$. For each pair of neighbours of $v$, there must be two vertices in $S$ which are common neighbours of both of them. Denote these subsets of $S$ by $C(w,x)$, $C(w,y)$, etc.

The sets $C(w,x)$ and $C(w,y)$ cannot be disjoint because $w$ has degree 4. If $C(w,x)=C(w,y)$, then since $w$ has degree 4, we must also have $C(w,z)\cap C(w,x)\neq \emptyset$, but then a vertex in $C(w,z)\cap C(w,x)$ has four common neighbours with $v$, a contradiction. 

So each pair of the six considered subsets of $S$ intersect in one vertex, and each of these vertices belongs to at most 3 of the subsets. This happens only if $|S|=4$ and each element of $S$ is adjacent to a distinct set of three vertices from $\{w,x,y,z\}$. Finally, each element of $S$ has one further neighbour, but each pair of vertices in $S$ has one further common neighbour, hence $S$ is the neighbourhood of a single vertex $v'$, and $G$ is isomorphic to $K_{5,5}-M$.
\end{proof}

\begin{theorem}
    If $D$ is a connected regular oriented graph with exactly two distinct $H_{\omega}$-eigenvalues $r>s$, then the following are equivalent:
    \begin{enumerate}
        \item $r=2$.
        \item $D$ has no underlying triangles.
        \item The underlying graph of $D$ is isomorphic to $K_{5,5}-M$.
    \end{enumerate}
\end{theorem}

\begin{proof}
$(3) \Rightarrow (1)$ and $(3)\Rightarrow(2)$ are trivial. $(2)\Rightarrow (1)$ is Lemma \ref{Lem:r=2}. If $(2)$ holds, then so does $(1)$, hence $r=-s=2$ and so the assumptions of Lemma \ref{Lem:3common} hold. Thus the assumptions of Lemma \ref{Lem:k55} hold, which implies $(3)$, so we have $(2)\Rightarrow (3)$.

Finally, suppose $(1)$ holds but $(2)$ does not. Then consider the bipartite double of $D$, that is, the oriented graph $\widetilde{D}$ with Hermitian adjacency matrix \[\pmat{0 &H_\omega(D)\\H_\omega(D) &0}.\] This is a triangle-free, regular oriented graph. Moreover, the underlying graph of $\widetilde{D}$ is connected since the underlying graph of $D$ is connected and not bipartite (see Theorem 3.4 from \cite{BHM80}). For any eigenvector $v$ of $H_\omega(D)$ with eigenvalue $\theta$, the vectors $(v,v)$ and $(v,-v)$ are eigenvectors of $H_\omega(\widetilde{D})$ with eigenvalues $\theta$ and $-\theta$ respectively. The set \[\{(v,v),(v,-
v): v \mbox{ is an eigenvector of } H_\omega(D)\}\] spans the eigenvectors of $H_\omega(\widetilde{D})$, hence $H_\omega(\widetilde{D})$ has exactly two distinct eigenvalues, $2$ and $-2$. Now, since $(2)\Rightarrow (3)$ and since $(2)$ holds for $\widetilde{D}$, the underlying graph of $\widetilde{D}$ must be isomorphic to $K_{5,5}-M$. This in turn implies that the underlying graph of $D$ is isomorphic to $K_5$. Up to isomorphism there is only a single regular tournament of order 5. We computed its $H_\omega$-eigenvalues directly and found that they are all distinct, a contradiction.
\end{proof}

For the case $r=2$, all that is left is to check the possible orientations of $K_{5,5}-M$. 

\begin{lemma}
    Up to isomorphism there is a unique orientation of $K_{5,5}-M$ which is regular and has two $H_\omega$-eigenvalues.
\end{lemma}

\begin{proof}
Label the vertices of one part $v_1,w_1,x_1,y_1,z_1$ and the vertices of the other part $v_2,w_2,x_2,y_2,z_2$ with $a_1$ non-adjacent to $a_2$ for $a\in\{v,w,x,y,z\}$. Without loss of generality we may assume $v_1\rightarrow w_2$, $v_1\rightarrow x_2$, $v_1\leftarrow y_2$, $v_1 \leftarrow z_2$. Now, for two vertices $a,b$ from the same part of the bipartition, let $K(a,b)$ denote the $K_{2,3}$ subgraph induced by $\{a,b\}\cup (N(a)\cap N(b))$. From the proof of Lemma \ref{Lem:3common} we know that each $K_{2,3}$ subgraph has one of only two possible orientation. Looking at the partially oriented subgraph $K(v_1,w_1)$ we see that the only way for there to be a directed 2-walk from $v_1$ to $w_1$ is if $x_2\rightarrow w_1$. Similarly, looking at the subgraphs $K(v_1,x_1)$, $K(v_1,y_1)$, $K(v_1,z_1)$ we find that $w_2\rightarrow x_1$, $z_2\leftarrow y_1$ and $y_2 \leftarrow z_1$ respectively. 

The remaining undirected edges consist of the four edges incident with $v_2$ and the edges of two 4-cycles: $(w_2,y_1,x_2,z_1)$ and $(w_1,y_2,x_1,z_2)$. Since $w_2$ and $x_2$ already have one in-neighbour and one out-neighbour the orientation on the 4-cycle $(w_2,y_1,x_2,z_1)$ gives either a directed cycle (in either direction) or it has two arcs into $y_1$ and two arcs out of $z_1$, or it has two arcs into $z_1$ and two arcs out of $y_1$. But $y_1$ and $z_1$ each already have one arc leaving, so the later two cases are impossible, hence $(w_2,y_1,x_2,z_1)$ must be a directed cycle, which implies $y_1\leftarrow v_2$ and $z_1\leftarrow v_2$. By a similar argument, $(w_1,y_2,x_1,z_2)$ is a directed cycle which implies $w_1\rightarrow v_2$, and $x_1 \rightarrow v_2$.

It remains to determine which direction each of the two 4-cycles go in, but these four choices yield four isomorphic oriented graphs: The map which exchanges the pairs $\{w_1,x_1\}$ and $\{y_2,z_2\}$ reverses the order of one of the 4-cycles. The map which exchanges the pairs $\{w_2,x_2\}$ and $\{y_1,z_1\}$ reverses the other 4-cycle.
\end{proof}

\subsection{Regular Oriented Graphs with $r\geq 3$} \label{Sec:rgeq3}

The remaining candidates for oriented graphs with two distinct $H_{\omega}$-eigenvalues are $r$-regular and have eigenvalues $r$ and $-2$. We show that such oriented graphs do not exist. This will complete the proof of Theorem \ref{Thm:2ev_Char}.

\begin{lemma}\label{Lem:rleq2}
    Let $D$ be an $r$-regular oriented graph with two distinct $H_{\omega}$-eigenvalues $r$ and $s$, where $r>s$. If $s=-2$, then $r\leq 2$.
\end{lemma}

\begin{proof}
    Let $H=H_\omega(D)$. Since $s=-2$, the minimal polynomial of $H$ is $(x-r)(x+2)$ and therefore we have 
    \begin{equation} \label{eq:Min_poly2}
        H^2-(r-2)H-2rI=0.
    \end{equation}
    Pick an arc $u\rightarrow v$. Examining the $(u,v)$ entry in the matrix equation (\ref{eq:Min_poly2}), we see that the sum of the values of all $2$-walks from $u$ to $v$ is $(r-2)\omega$.
    
    The value of any $2$-walk is in $\{1,\omega^2,\omega^4\}$, so we have 
    \[a+b\omega^2+c\omega^4=(r-2)\omega\] 
    where $a,b,c$ are non-negative integers counting the number of 2-walks from $u$ to $v$ with values $1$, $\omega^2$ and $\omega^4$ respectively. Subtracting $a(1+\omega^2+\omega^4)=0$ from the above equation, multiplying both sides by $\omega^{-1}$ and rearranging we have 
    \[(b-a)\omega=(a-c)\omega^3+r-2.\] 
    Since $\omega^3=-1$ the right-hand side is real, so the left-hand side is zero. We conclude that $a=b$ and $a-c=r-2.$ Since $c$ is non-negative, $a\geq r-2$, so there are at least $r-2$ walks of length $2$ with value $1$ from $u$ to $v$. Since $b=a$ there are also at least $r-2$ walks of length $2$ with value $\omega^2$ from $u$ to $v$. 
    
    To show that $r\leq4$ we argue that there are at most two $2$-walks of value $1$ from $u$ to $v$.
    The orientations of $2$-walks from $u$ to $v$ with value $1$ are $u\rightarrow x\leftarrow v$ and  $u\leftarrow x\rightarrow v$, which we will call \textbf{absorbing} and \textbf{repelling}, respectively.
    Since the only $2$-walk with value $\omega^2$ is the directed $2$-walk, there must be $r-2$ directed $2$-walks from $u$ to $v$.
    In particular, these paths contribute $r-2$ incoming edges to $v$, and combined with the arc $(u,v)$, this means $v$ can only have one further incoming arc.
    Thus, there must be at most one repelling $2$-walk from $u$ to $v$.
    Similarly, the $r-2$ directed 2-walks contribute $r-2$ outgoing edges from $u$, and combined with the arc $(u,v)$, this means there must be at most one absorbing 2-walk from $u$ to $v$. 
    
    It remains to show that $r=3$ and $r=4$ are not possible.

    \textbf{Case $r=4$:} As above, pick an arc $u\rightarrow v$. Since $r=4$ there are two directed 2-paths from $u$ to $v$: $u\rightarrow x_1\rightarrow v$ and $u \rightarrow x_2 \rightarrow v$. There is also one absorbing 2-walk and one repelling 2-walk: $u\rightarrow y_1 \leftarrow v$ and $u \leftarrow y_2 \rightarrow v$. The nine arcs just listed include for all four of the arcs into $v$. 

    Since $x_1\rightarrow v$, there must be two directed 2-paths from $x_1$ to $v$, and since each arc into $v$ is accounted for, this implies $x_1\rightarrow x_2$ and $x_1\rightarrow y_2$. But since $x_2\rightarrow v$, we must also have $x_2\rightarrow x_1$ by the same argument. But the edge $\{x_1,x_2\}$ cannot be oriented in both directions simultaneously, a contradiction. 
    
    \textbf{Case $r=3$:} Since $r=3$ there is at least one directed $2$-walk from $u$ to $v$, so there is a vertex $x$ such that $u\rightarrow x \rightarrow v$. Moreover there is a vertex $y$ adjacent to both $u$ and $v$ so that the 2-walk $u,y,v$ has value $1$. 
    
    \textbf{Subcase $1$:} The walk $u,y,v$ is absorbing: $u\rightarrow y \leftarrow v$. Since $u\rightarrow x$ and all of $u$'s outgoing edges are accounted for, we must have $y\rightarrow x$. Now we consider the 2-walks from $y$ to $x$. We already have $y\leftarrow u\rightarrow x$ of value 1 and $y\leftarrow v \leftarrow x$ of value $\omega^4$. Since $y\rightarrow x$ we must now have at least two additional directed 2-walks from $y$ to $x$, one to cancel with the existing 2-walks of values 1 and $\omega^4$, and one to contribute to $H_{x,y}=\omega$. But then $x$ must have in-degree at least 4, a contradiction.

    \textbf{Subcase $2$:} The walk $u,y,v$ is repelling: $u\leftarrow y \rightarrow v$. Since $x\rightarrow v$ and all of $v$'s incoming edges are accounted for, we must have $x\rightarrow y$. Now we consider the 2-walks from $x$ to $y$. We already have $x\rightarrow v \leftarrow y$ of value 1 and $x\leftarrow u \leftarrow y$ of value $\omega^4$. Since $x\rightarrow y$ we must now have at least two additional directed 2-walks from $x$ to $y$, one to cancel with the existing 2-walks of values 1 and $\omega^4$, and one to contribute to  $H_{y,x}=\omega$. But then $x$ must have out-degree at least 4, a contradiction.
    \end{proof}

\subsection{Extension to Mixed Graphs}

In this section we sketch the extension of Theorem \ref{Thm:2ev_Char} to mixed graphs. The arguments are a simple extension of the above and yield only one additional non-trivial example.

Recall that a \textbf{mixed graph} $D=(V,A,E)$ consists of a vertex set $V$ a set $A$ of arcs and a set $E$ of edges. As described in the introduction we extend the definition of $H_\omega(D)$ so that $H_\omega(D)_{u,v}=1$ when $\{u,v\}$ is an edge of $D$.

\begin{theorem} \label{Thm:Mixed_2ev}
If $D$ is a connected mixed graph and $H_\omega(D)$ has exactly two distinct eigenvalues, then $D$ is one of the following.

\begin{enumerate}
    \item[(a)] An complete graph on $n\geq 2$ vertices cospectral to the unoriented complete graph.
    \item[(b)] A 4-cycle containing a directed path of length three and one edge.
    \item[(c)] One of the four oriented graphs from Theorem \ref{Thm:2ev_Char}. 
    
\end{enumerate}
\end{theorem}

\begin{proof}
As before we let $H=H_\omega(D)$ and let $s<0<r$ be the two distinct eigenvalues of $H$.

Note that Lemma \ref{Lem:min_poly} holds for mixed graphs. The only difference in the proof is when showing that $(r+s)$ must be an integer. There we must consider five types of walks of length 2 between $u$ and $v$ instead of three, but the same technique applies.

Next, note that the bound in Lemma \ref{Lem:s_geq2} still holds and is sharp if $D$ has any undirected edges. So if $D$ is a mixed graph we have $s\in \{-1,-\sqrt{2},-\sqrt{3}\}$. 

Suppose $s=-1$. Interlacing implies that each three vertex induced subgraph of $D$ is a triangle with at least one edge and no directed path of length 2. It follows that we are in case $(a)$.

Suppose $s=-\sqrt{2}$. The graph underlying $D$ is cycle.
Lemma \ref{Lem:3common} no longer holds but the argument in its proof implies that any pair of vertices at distance 2 in $D$ must have a number of common neighbours equal to $2a+3b$ for non-negative integers $a$ and $b$. So $C_4$ is the only new possibility for the underlying graph. The mixed orientation stated in case $(b)$ is easily verified to be the unique choice with two eigenvalues. 

Suppose $s=-\sqrt{3}$. Note first that each mixed orientation of the triangle has an eigenvalue outside the interval $[-\sqrt{3},\sqrt{3}]$, so interlacing implies the graph underlying $D$ is triangle-free. Now in the underlying graph all vertices have degree 3 and each pair of vertices at distance 2 must have two or three common neighbours. 

If there is a pair of vertices with three common neighbours then the argument used in the proof of Theorem \ref{Thm:K33} shows that $D$ must be the oriented $K_{3,3}$ found there. So we assume all vertices of $D$ have two common neighbors. 


It is easy to verify that the cube is the only connected cubic triangle-free graph in which each pair of vertices at distance 2 have two common neighbors. 

Suppose we have a mixed orientation of the cube with two $H_\omega$-eigenvalues. Each of the six 4-cycles must be the mixed orientation of $C_4$ mentioned above which contains one undirected edge. So this mixed cube has exactly three undirected edges and they form a maximal matching. The two remaining vertices uncovered by this matching are each required to be the center of three directed paths of length 2, one for each incident 4-cycle, but this is impossible. \end{proof}

\section{Three Distinct Eigenvalues}\label{Sec:3ev}
In counterpoint to the previous section, we show that oriented graphs with three distinct $H_{\omega}$-eigenvalues are plentiful. The results of this section are a generalization of results from \cite{dCGKPM92}. 

First, note that since $\omega=\frac{1+i\sqrt{3}}{2}$ and $\overline{\omega}=\frac{1-i\sqrt{3}}{2}$ the Hermitian adjacency matrix of any tournament $T$ is \[H_\omega(T)=\frac{1}{2}(J-I)+i\frac{\sqrt{3}}{2}(B-I)\] where $J$ is the all ones matrix and $B-I$ is a skew-symmetric $(1,-1)$-matrix. 

A \textbf{Hadamard matrix} is a square $(1,-1)$-matrix of order $n$ such that $HH^T=nI$. A \textbf{skew-Hadamard matrix} is a Hadamard matrix for which $H-I$ is skew-symmetric. Many constructions of skew-Hadamard matrices are known, see \cite{KS08} for a survey of constructions. For us, the important fact is that there are infinitely many integers $k$ for which a skew-Hadamard matrix of order $4k$ exists. A \textbf{regular tournament} is a regular oriented complete graph.

\begin{theorem}
For each positive integer $k\geq2$ for which there is a skew-Hadamard matrix of order $4k$, there is a regular tournament of order $4k-1$ for which $H_\omega(T)$ has exactly three distinct eigenvalues.
\end{theorem}
\begin{proof}
    Let $A$ be a skew-Hadamard matrix of order $n=4k$. Without loss of
    generality we may assume that the first row of $A$ is $j^\intercal$, the all one's vector. Now, remove the first row and the first column of $A$ to obtain a matrix $B$ so that  $B-I$ is a skew-symmetric matrix of order $4k-1.$ Note that in every row and column  of $B-I$ there are $2k-1$ entries  1 and $2k-1$ entries $-1$. Clearly, $BB^\intercal=nI-J$ and $BJ=JB$. Since $B-I$ is skew-symmetric, we have $(B-I)^\intercal=I-B$ and so we find that $B+B^\intercal=2I$.  
    
    Now let $T$ be the $(2k-1)$-regular tournament whose Hermitian adjacency matrix is  \[H=H_\omega(T)=\frac{J-I}{2}+i\frac{\sqrt{3}}{2}(B-I).\] Note that $T$ has order $4k-1$. We have \[(B-I)^2=B^2-2B+I=B(B-2I)+I=-BB^\intercal+I=J-(n-1)I.\]

     It follows that the all one's vector of length $n-1$ is an eigenvector of $B-I$ with eigenvalue 0, and the remaining eigenvectors are orthogonal to $J$, hence the remaining eigenvalues of $B-I$ are $\pm i \sqrt{n-1}$. Since $\tr(B-I)=0$, this latter pair of eigenvalues have equal multiplicity $\frac{n-2}{2}$.

     Since $BJ=JB$, the three matrices $B-I, J$ and $I$ are simultaneously diagonalizable (note that $B$ is diagonalizable because $iB$ is Hermitian). Thus the eigenvalues of $H$ are $\frac{n-2}{2}$ of multiplicity 1 and \[-\frac{1}{2}\pm\frac{\sqrt{3(n-1)}}{2},\] each with multiplicity $\frac{n-2}{2}$. Hence $T$ is of order $n-1$ and
    has 3 distinct eigenvalues and the proof is complete.
    \end{proof}

In the other direction, any regular tournament with three $H_{\omega}$-eigenvalues produces a skew-Hadamard matrix as well.

\begin{theorem}
Let $T$ be a regular tournament of order $n-1$ with $3$ distinct $H_{\omega}$-eigenvalues. Then there is a skew-Hadamard matrix of order $n$.
\end{theorem}

\begin{proof}
Let $H$ be the adjacency matrix of $T$. Since $T$ is regular, we have
\[H=\frac{J-I}{2}+i\frac{\sqrt{3}}{2}(B-I),\] where $B-I$ is a skew-symmetric $(1,-1)$-matrix, whose all diagonal elements are zero. Since $T$ is a regular tournament of order $n-1$, one of the eigenvalues of $T$ is $\frac{n-2}{2}$.
If $\theta$ is an eigenvalue of $B-I$, then $-\theta$
is also an eigenvalue of $B-I$. Thus the three distinct eigenvalues of $H$ are $\frac{n-2}{2}$, with multiplicity 1, and $-\frac{1}{2}\pm i\frac{\sqrt{3}}{2} \theta$, each with multiplicity $\frac{n-2}{2}$. Since $\tr(H^2)$ is twice the number of edges we have 
\[\Big(\frac{n-2}{2}\Big)^2+ \frac{n-2}{2} \Big(-\frac{1}{2}+ i\frac{\sqrt{3}}{2} \theta\Big)^2+\frac{n-2}{2} \Big(-\frac{1}{2}-  i\frac{\sqrt{3}}{2} \theta\Big)^2=(n-1)(n-2),\] 
hence we find that $\theta=i\sqrt{n-1}$. 

Since $i(B-I)$ is a Hermitian matrix, it is diagonalizable and so
$B-I$ is diagonalizable. Thus the minimal polynomial of $B-I$ is
$x(x-\sqrt{n-1}\, i)(x+\sqrt{n-1}\, i)$.
This implies that \[(B-I)((B-I)^2+(n-1)I)=0.\] Zero is an eigenvalue of multiplicity 1 for the matrix $B-I$ and $j$ is its corresponding eigenvector. Moreover $(B-I)^2+(n-1)I$ is a symmetric matrix such that
its every column is an eigenvector of zero for the matrix $B-I$, so we conclude that $(B-I)^2+(n-1)I=cJ$, for some suitable $c$.

 Since $(B-I)j=0$, we conclude that $(n-1)j=cJj=c(n-1)j$. Therefore $c=1$ and  $(B-I)^2=J-(n-1)I$. Since $B-I$ is a skew-symmetric matrix, we have \[(B-I)(B-I)^\intercal=-(B-I)^2=(n-1)I-J.\] This equation implies that the matrix $A$ obtained from $B$ by adding a column of $-1$'s  and then a row of $1$'s is a skew-Hadamard matrix of order $n$. 
 \end{proof}

It would be interesting to know if the assumption that $T$ is regular can be removed.

\section{Other Roots of Unity}\label{Sec:other_roots}
Let $\sigma$ be a primitive $k$-th root of unity for some positive integer $k$. We consider here the possibility of characterizing the oriented graphs $D$ so that $H_\sigma(D)$ has exactly two distinct eigenvalues.

\subsection{$k=3$}
Let $\gamma$ be a primitive third root of unity. Since the matrices $H_\gamma(D)$ and $H^T_\gamma(D)$ are similar we may assume, without loss of generality, that $\gamma=-\frac{1}{2}-i\frac{\sqrt{3}}{2}$. This implies that, for any oriented graph, $H_\gamma(D)=-H_\omega(D)$. So for oriented graphs the $k=3$ variant of our problem is essentially identical to the $k=6$ case that we characterized in Section \ref{Sec:two_ev}.

\begin{corollary}
      If $D$ is a connected oriented graph and $H_\gamma(D)$ has exactly two distinct eigenvalues, then $D$ is one of the four oriented graphs from Theorem \ref{Thm:2ev_Char}.
\end{corollary}

The results of Section \ref{Sec:3ev} extend to $H_\gamma$ just as easily.

However, for mixed graphs the situation is very different. For each even $n$ there are non-trivial example of mixed orientations of $K_{3^n}$ with exactly two $H_\gamma$-eigenvalues. See Section 12.5 of \cite{BCN89}.

\subsection{$k=4$}
This case bares great similarity to the case of signed graphs. In fact, we next show that a large class of signed graphs with two eigenvalues is equivalent to oriented graphs with the same two $H_i$ eigenvalues.

Recall that a \textbf{signed graph} is a pair $(G,\phi)$ consisting of graph $G$ and a map $\phi:E(G)\rightarrow\{1,-1\}.$ The adjacency matrix $S=S(G)$ of a signed graph is defined by \[S_{u,v}=\begin{cases}
    1&   \phi(\{u,v\})=1 \\
    -1&\phi(\{u,v\})=-1 \\
    0&\mbox{otherwise.}
\end{cases}\]

\begin{theorem} \label{Thm:signed_to_i}
    Let $G$ be a simple undirected bipartite graph. The following are equivalent.
    \begin{itemize}
        \item There exists a signing $\phi$ so that the adjacency matrix of the signed graph $(G,\phi)$ has exactly two distinct eigenvalues.
        \item There is an orientation $D$ of $G$ so that $H=H_i(D)$ has exactly two distinct eigenvalues.

    \end{itemize}
\end{theorem}

\begin{proof}
    First suppose there is such a signed graph with signed adjacency matrix \[S=\pmat{0&B\\B^T&0}.\] 
    
    We calculate that \[\pmat{iI & 0 \\0 &I}\pmat{0&B\\B^T&0}\pmat{-iI&0\\0&I}=i\pmat{0&B\\-B^T&0}=H,\] where $H=H_i(D)$ is the Hermitian adjacency matrix of some oriented bipartite graph $D$. So $S$ and $H$ are similar, and $D$ has two $H_i$-eigenvalues. The converse direction is very similar.
\end{proof}

Theorem \ref{Thm:signed_to_i} applies to many known signed graphs. Most notably, to an infinite family of 4-regular signed graphs (see \cite{HTW19}), and an infinite family of signed hypercubes, see \cite{HH19}.

It is possible that the assumption that $G$ is bipartite can be relaxed somewhat, but it cannot be dropped altogether: There is a signed $K_6$ with eigenvalues $\pm\sqrt{5}$, but we checked, with a computer, all orientations of $K_6$ and found that none of them yield a matrix $H_i(D)$ with two eigenvalues.

\subsection{Other values of $k$}

Now let $\sigma$ be any primitive $k$th root of unity. Certain aspects of our work in Section 
\ref{Sec:two_ev} remain applicable: Our first key result, Lemma \ref{Lem:s_geq2}, carries over nicely as long as $\mbox{Re}(\sigma)>0$.

\begin{corollary}\label{Cor:s>re}
 Let $D$ be an oriented graph with two distinct $H_{\sigma}$-eigenvalues $r$ and $s$, with $r>s$. Suppose $\mbox{Re}(\sigma)>0$. Then $s\geq\frac{-1}{\mbox{Re}(\sigma)}$, with equality if and only if $D$ is $r$-regular.
\end{corollary}

\begin{proof}
    The minimal polynomial of $H$ is still quadratic, so equations (\ref{Eq:Min_poly1})  
    and (\ref{eq:nrs}) from the proof of Lemma \ref{Lem:s_geq2} remain valid. Now we have \[r\geq \frac{j^\intercal Hj}{j^\intercal j}=\frac{1}{n}\sum_{i>j}(H_{i,j}+H_{j,i})=\frac{2|E(D)|\mbox{Re}(\sigma)}{n}=\frac{(-nrs)\mbox{Re}(\sigma)}{n}\] with equality if and only if $D$ is regular. The bound on $s$ follows.
\end{proof}

If we restrict to the primitive $k$th root of unity $\sigma_1$ with maximal real part, then $\mbox{Re}({\sigma_1})$ is an increasing function of $k$ and the above bound becomes very strong.

\begin{theorem}
    Let $k>8$ be an integer and let $\sigma_1=\cos(\frac{2\pi}{k})+i\sin(\frac{2\pi}{k})$. The orientation of $K_2$ is the only connected oriented graph $D$ such that $H_{\sigma_1}(D)$ has exactly two distinct eigenvalues.
\end{theorem}

\begin{proof}
    First suppose $D$ is not a tournament. Then it contains an induced path on three vertices. Any orientation of the path on three vertices has $H_{\sigma_1}$-eigenvalues $\{-\sqrt{2},0,\sqrt{2}\}$, regardless of the value of $k$. So interlacing implies that the least eigenvalue of $H_{\sigma_1}(D)$ is at most $-\sqrt{2}$. But $\mbox{Re}({\sigma_1})=\cos(\frac{2\pi}{k})$ is an increasing function of $k$ whose value at $k=8$ is $\frac{1}{\sqrt{2}}$. So for any larger $k$, the bound from Corollary \ref{Cor:s>re} implies $s>-\sqrt{2}$, a contradiction.

    Now suppose $D$ is a tournament. Let $r>s$ denote the two distinct eigenvalues. Suppose the multiplicity of $r$ is greater than one.  Corollary \ref{Cor:s>re} implies $r\geq 2\cos(\frac{2\pi}{k})d=(n-1)\cos(\frac{2\pi}{k})$ and $s\geq -\frac{1}{\cos(\frac{2\pi}{k})}$. So we have
    
    \[0=\tr(H_{\sigma_1}(D))=\sum_{i=1}^n \lambda_i\geq 2(n-1)\cos(\frac{2\pi}{k})-(n-2)\frac{1}{\cos(\frac{2\pi}{k})}.\] But since $k>8$, $\cos(\frac{2\pi}{k})>\frac{1}{\sqrt2}$, so  \[0\geq \frac{2(n-1)}{(n-2)}\cos^2(\frac{2\pi}{k})-1\geq\frac{n-1}{n-2}-1>0,\] a contradiction.

    So $r$ must have multiplicity 1 which implies $s$ has multiplicity $n-1$. The $H_\sigma$-adjacency matrix of an oriented $K_2$ has eigenvalues $\{1,-1\}$ and these eigenvalues must interlace the eigenvalues of $H_{\sigma_1}(D)$, but since $s$ has multiplicity $n-1$ this implies $s\geq -1 \geq s$, so $s=-1$. But the $H_{\sigma_1}$-adjacency matrices of the two distinct orientations of the triangle each have an eigenvalue less than $-1$, so there are no induced triangles in our tournament and $n=2$. \end{proof}

In the opposite extreme, for a primitive $k$th root $\sigma_2$ with $\mbox{Re}(\sigma_2)$ close to $-1$, one can use the fact that \[s\leq \frac{j^\intercal Hj}{j^\intercal j}=2\mbox{Re}(\sigma_2)d\] to get a strong upper bound on $r$ when $k$ is large.

The intermediary choices, where $\mbox{Re}(\sigma)$ is close to zero, yield weaker bounds. This is most evident with $k=4$ where there exist oriented $n$-cubes which have $H_i$-eigenvalues $\pm \sqrt{n}$. When $\mbox{Re}(\sigma)$ is near zero we seem to have the best chance of finding large families of examples, but we do not yet know if there are infinite families when $k\neq 4$. So we conclude with a question.

Is there a primitive root of unity $\sigma$ with $\mbox{Re}(\sigma)\neq 0$ for which there exist infinitely many oriented graphs $D$ with exactly two distinct $H_\sigma(D)$-eigenvalues?

\section*{Conflict of Interest}

Given his role as Associate Editor, Bojan Mohar had no involvement in the peer review of this article and had no access to information regarding its peer review. Full responsibility for the editorial process for this article was delegated to another journal editor.

\bibliographystyle{abbrv}
\bibliography{HermMat}

\end{document}